\documentclass[twoside,12pt]{article}
\usepackage[]{amsmath,amsthm,amssymb}
\usepackage[latin1]{inputenc}

\begin{document}

\title{{\Large\bf  Discrete Mehler-Fock transforms}}

\author{Semyon  YAKUBOVICH}
\maketitle

\markboth{\rm \centerline{ Semyon   YAKUBOVICH}}{}
\markright{\rm \centerline{MEHLER-FOCK  TRANSFORMS}}

\begin{abstract} {\noindent Discrete analogs of the classical Mehler-Fock  transforms are introduced and investigated. It involves series with the associated Legendre  function $P^\mu_{in-1/2}(x), x > 1,\  {\rm Re} \mu < 1/2, \ n \in \mathbb{N}, i $ is the imaginary unit, and the so-called incomplete Legendre  integrals. Several expansions of arbitrary functions and sequences in terms of these series and integrals are established.  }

\end{abstract}
\vspace{4mm}



\vspace{4mm}

\section {Introduction and preliminary results}

In 1943  V.A. Fock   (cf. [2]) established   the following expansion of an arbitrary function $f$ in terms of the spherical Legendre functions $P_{i\tau-1/2} (x),\ \tau > 0,\ x  > 1$ [1], Vol. I

$$f(\tau) =  \tau \tanh(\pi\tau) \int_1^\infty P_{i\tau-1/2} (x) \int_0^\infty P_{iy-1/2}(x) f(y)\ dy dx,\ \tau >0,\eqno(1.1)$$
which generates classical  Mehler-Fock transforms [5].  Later in 1958 N.Ya.Vilenkin [4] generalized (1.1) to the associated Legendre functions $P^\mu_{i\tau-1/2} (x),\mu \in \mathbb{C},\   P^0_{i\tau-1/2} (x) \equiv P_{i\tau-1/2} (x)$  [1], Vol. I as the kernel

$$f(\tau) =  {\tau \over \pi}  \sinh(\pi\tau) \Gamma\left({1\over 2}+ i\tau-\mu\right) \Gamma\left({1\over 2}- i\tau-\mu\right)$$

$$\times \int_1^\infty P^\mu_{i\tau-1/2} (x) \int_0^\infty P^\mu_{iy-1/2}(x) f(y)\ dy dx,\ \tau >0,\eqno(1.2)  $$
where $\Gamma(z)$ is Euler's gamma-function  [1], Vol. I.  The associated Legendre functions $P^\mu_{\nu} (z)$ are solutions of the corresponding second order differential Legendre equation with polynomial coefficients

$$(1-z^2) {d^2\omega\over dz^2} - 2z  {d\omega\over dz} + \left(\nu(1+\nu)-\mu^2(1-z^2)^{-1}\right)\omega= 0,\eqno(1.3)$$
and can be expressed in terms of the Gauss hypergeometric function [1], Vol. I

$$\omega= P^\mu_{\nu} (z) = {1\over \Gamma(1-\mu)} \left({z+1\over z-1}\right)^{\mu/2} {}_2F_1\left(-\nu,\ 1+\nu;\ 1-\mu;\ {1-z\over 2}\right), \quad |1-z| < 2.$$
Besides, it can be represented by the Mehler and Legendre integrals, respectively,

$$P^\mu_{\nu} (\cosh(\alpha)) = \sqrt{{2\over \pi}} { (\sinh(\alpha))^\mu\over \Gamma(1/2-\mu)} \int_0^\alpha {\cosh(t(\nu+1/2))\over (\cosh(\alpha) - \cosh(t))^{1/2+\mu}}\ dt,\eqno(1.4)$$
where $\alpha >0,\  {\rm Re} (\mu) < 1/2$, 

$$P^\mu_{\nu} (z) = \sqrt{{2\over \pi}}\  {\Gamma(1/2-\mu) (z^2-1)^{-\mu/2} \over \Gamma(1+\nu-\mu) \Gamma(-\nu-\mu)} \int_0^\infty {\cosh(t(\nu+1/2))\over (\cosh(t) + z)^{1/2-\mu}}\ dt,\ z >1\eqno(1.5)$$
under conditions ${\rm Re} (\nu+\mu) < 0,\  {\rm Re} (\nu-\mu+1) > 0$, and in terms of the Mellin-Barnes integral (see [3], Vol. III, Entry 8.4.41.10), which will be used in the sequel

$$  \Gamma(1+\nu-\mu) \Gamma(-\nu-\mu) (1+x)^{\mu/2} P^\mu_{\nu} \left(2x+1\right)$$

$$ = {1\over 2\pi i} \int_{\gamma-i\infty}^{\gamma-i\infty} \frac{ \Gamma(s-\mu/2) \Gamma(1+\nu-\mu/2-s)
\Gamma(-\nu-\mu/2-s)}{\Gamma(1-\mu/2-s)} x^{-s} ds,\eqno(1.6)$$
where $x >0,\ {\rm Re} \mu/2 < \gamma < 1+ {\rm Re} (\nu-\mu/2),  \ - {\rm Re} (\nu+\mu/2).$

  The main goal of this paper is to establish the following expansions of arbitrary  sequences and functions which give rise to  discrete analogs of the Mehler-Fock transforms (1.1), (1.2), namely

$$a_n =   {n\over \pi}  \sinh(\pi n) \Gamma\left({1\over 2} + in -\mu\right) \Gamma\left({1\over 2} - in -\mu\right) $$

$$\times \int_1^\infty   P^\mu_{in-1/2}(x, \pi) \sum_{m=1}^\infty a_m  P^\mu_{im-1/2}(x) dx,\ n \in \mathbb{N},\eqno(1.7)$$

$$a_n =   {n\over \pi}  \sinh(\pi n) \int_1^\infty    P^\mu_{in-1/2}(x)  $$

$$\times \sum_{m=1}^\infty a_m  \Gamma\left({1\over 2} + im -\mu\right) \Gamma\left({1\over 2} - im -\mu\right) P^\mu_{im-1/2}(x, \pi)  dx,\  n \in \mathbb{N},\eqno(1.8)$$

$$ f (x) =   {1\over \pi}     \sum_{n=1}^\infty  n\sinh(\pi n) \Gamma\left({1\over 2} + in -\mu\right) \Gamma\left({1\over 2} - in -\mu\right) P^\mu_{in-1/2}(x, \pi)  $$

$$\times \int_1^\infty  P^\mu_{in-1/2}(t) f(t) dt,\ x > 1,\eqno(1.9)$$

$$ f (x) =   {1\over \pi}     \sum_{n=1}^\infty  n\sinh(\pi n) \Gamma\left({1\over 2} + in -\mu\right) \Gamma\left({1\over 2} - in -\mu\right) P^\mu_{in-1/2}(x)  $$

$$\times \int_1^\infty  P^\mu_{in-1/2}(t,\pi) f(t) dt,\ x > 1.\eqno(1.10)$$
Here the function  $P^\mu_{\nu}(z, \omega),\ z >1,\ \omega >0$ denotes the incomplete Legendre integral 

$$P^\mu_{\nu}(z, \omega) =  \sqrt{{2\over \pi}}\  {\Gamma(1/2-\mu) (z^2-1)^{-\mu/2} \over \Gamma(1+\nu-\mu) \Gamma(-\nu-\mu)} \int_0^\omega  {\cosh(t(\nu+1/2))\over (\cosh(t) + z)^{1/2-\mu}}\ dt.\eqno(1.11)$$
Finally in this section we note that the discrete analogs of the Kontorovich-Lebedev transform were considered recently by the author in [6]. It involves the modified Bessel functions and incomplete Bessel integrals [1], Vol. II.  As is known [5], the Kontorovich-Lebedev and Mehler-Fock transforms   generate the theory of the index transforms, which is related to the integration with respect to parameters of the hypergeometric functions. 

\section{Expansion theorems}

We begin with 

{\bf Theorem 1.} {\it Let ${\rm Re} \mu < 1/2$ and the sequence $a= \{a_n\}_{n\in \mathbb{N}} \in l_1$, i.e. the following series converges

$$||a||= \sum_{n=1}^\infty  |a_n|   < \infty.\eqno(2.1)$$
Then  expansion $(1.7)$ holds, where the iterated series and integral converge absolutely.}

\begin{proof}  Indeed, denoting by 

$$F(x)=  \sum_{m=1}^\infty a_m  P^\mu_{im-1/2}(x),\eqno(2.2)$$
we see that the series converges absolutely for any $x >1$ by virtue of the estimate (see (1.4))

$$ |F(x)| \le  \sum_{m=1}^\infty  \left| a_m  P^\mu_{im-1/2}(x)\right| \le  \left| P^\mu_{-1/2}(x)\right| \sum_{m=1}^\infty  \left| a_m \right|.\eqno(2.3)$$
Hence, multiplying both sides of (2.2) by $(x^2-1)^{-\mu/2} e^{-yx},\ y >0$, we integrate over $(1,\infty)$ to get 

$$\int_1^\infty  (x^2-1)^{-\mu/2} e^{-yx} F(x) dx =  \sum_{m=1}^\infty a_m  \int_1^\infty  (x^2-1)^{-\mu/2} e^{-yx}  P^\mu_{im-1/2}(x) dx,\eqno(2.4)$$
where the interchange of the order of integration and summation is permitted by Fubini's theorem via the estimate

$$\int_1^\infty \left| (x^2-1)^{-\mu/2}\right|  e^{-yx} \sum_{m=1}^\infty \left| a_m  P^\mu_{im-1/2}(x)\right| dx \le $$

$$\le \int_1^\infty   e^{-yx} (x^2-1)^{-{\rm Re} \mu/2} \left| P^\mu_{-1/2}(x)\right| dx \sum_{m=1}^\infty \left| a_m  \right| < \infty$$
since $\left| P^\mu_{-1/2}(x)\right|$ can be  estimated, recalling integral (1.6). Namely, letting $x= 2t+1, t >0$, we have

$$ \left| P^\mu_{-1/2}(2t+1)\right| \le C_\mu {t^{-\gamma}\over (1+t)^{{\rm Re} \mu/2}},\   {{\rm Re} \mu\over 2} < \gamma < {1- {\rm Re} \mu\over 2},\eqno(2.5)$$
where 

$$C_\mu=  {1\over 2\pi} \int_{\gamma-i\infty}^{\gamma-i\infty} \left| \frac{ \Gamma(s-\mu/2) \Gamma^2 ((1-\mu)/2-s)
}{\Gamma^2(1/2-\mu)  \Gamma(1-\mu/2-s)}  ds \right|$$
and for each $ y >0$ the integral 

$$ \int_1^\infty   e^{-yx} (x^2-1)^{-{\rm Re} \mu/2} \left| P^\mu_{-1/2}(x)\right| dx = 2^{1- {\rm Re} \mu}  e^{-y} $$

$$\times  \int_0^\infty   e^{-yt} (t (t + 1))^{-{\rm Re} \mu/2} \left| P^\mu_{-1/2}(2t+1)\right| dt$$ 

$$\le  2^{1- {\rm Re} \mu}   C_\mu \int_0^\infty   e^{-yt} t^{-\gamma- {\rm Re} \mu/2} (t + 1)^{-{\rm Re} \mu} dt < \infty.$$
Returning to (2.4) and appealing to the formula in [3], Vol. III, Entry 2.17.7.5, we find

$$ \int_1^\infty  (x^2-1)^{-\mu/2} e^{-yx}  P^\mu_{im-1/2}(x) dx = \sqrt{{2\over \pi}}\  y^{\mu-1/2} K_{im}(y),\eqno(2.6)$$
where $K_\nu(z)$ is the modified Bessel function of the second kind [1], Vol. 2. This means that the right-hand side of (2.4) is the 
discrete Kontorovich-Lebedev transform [6] and the equality reads

$$ \sqrt{{\pi\over 2}} \  y^{1/2- \mu} \int_1^\infty  (x^2-1)^{-\mu/2} e^{-yx} F(x) dx =  \sum_{m=1}^\infty a_m  K_{im}(y).\eqno(2.7)$$
Moreover, the sequence $\{a_n\}_{n \ge 1}$ satisfies the condition of Theorem 1 in [6] and the inversion formula for the discrete Kontorovich-Lebedev transform holds.  This means that

$$a_n=   {\sqrt 2\over \pi\sqrt \pi} \ n  \sinh(\pi n) \int_0^\infty K_{in} (x,\pi)   x^{-1/2- \mu} \int_1^\infty  (t^2-1)^{-\mu/2} e^{-xt} F(t) dt,\eqno(2.8)$$
where $K_{in} (x,\pi)$ is the incomplete Bessel integral

$$K_{in} (x,\pi) = \int_0^\pi  e^{-x\cosh(u)} \cos (n u) du.\eqno(2.9)$$
But the estimates (2.3), (2.4) and Fubini's theorem justify the interchange of the order of integration in (2.8), and plugging  the integral (2.9) into (2.8), we derive with the use of (1.11)

$$ {\sqrt 2\over \pi\sqrt \pi} \ n  \sinh(\pi n) \int_0^\infty K_{in} (x,\pi)   x^{-1/2- \mu} \int_1^\infty  (t^2-1)^{-\mu/2} e^{-xt} F(t) dt$$

$$=  {\sqrt 2 \Gamma(1/2-\mu) \over \pi\sqrt \pi} \ n  \sinh(\pi n) \int_1^\infty  (t^2-1)^{-\mu/2} F(t) \int_0^\pi { \cos (n u) \over (t +\cosh(u))^{1/2-\mu}} du dt $$

$$=  {n \over \pi} \  \sinh(\pi n)  \Gamma\left({1\over 2} + in -\mu\right) \Gamma\left({1\over 2} - in -\mu\right) \int_1^\infty  P^\mu_{in-1/2}(t, \pi) F(t) dt ,$$
which proves (1.7).

\end{proof}

{\bf Remark 1}.  Expansion (1.7) generates a reciprocal pair of the generalized discrete Mehler-Fock transforms

$$F(x)=   \sum_{m=1}^\infty a_m  P^\mu_{im-1/2}(x),\ x >1,\eqno(2.10)$$

$$a_n=   {n\over \pi}  \sinh(\pi n) \Gamma\left({1\over 2} + in -\mu\right) \Gamma\left({1\over 2} - in -\mu\right) $$

$$\times \int_1^\infty   P^\mu_{in-1/2}(x, \pi) F(x) dx,\ n \in \mathbb{N}.\eqno(2.11)$$
 
{\bf Corollary 1}. {\it Let $\mu=0$ and condition $(2.1)$ hold. Then the following discrete analog of the classical Mehler-Fock expansion (1.1) takes place}

$$ a_n=   n  \tanh(\pi n) \int_1^\infty   P_{in-1/2}(x, \pi)   \sum_{m=1}^\infty a_m  P_{im-1/2}(x) dx,\ n \in \mathbb{N}.\eqno(2.12)$$

In order to establish expansion (1.8) we first will prove the following lemma about the orthogonality of the sequences of functions $\{ P^\mu_{in-1/2}(x)\}_{n \ge 1},\  \{  P^\mu_{in-1/2}(x,\pi)\}_{n \ge 1}$.  Precisely, it has 

{\bf Lemma 1.} {\it  Let $\left| {\rm Re} \mu \right| < 1/2$.  Then

$$\int_1^\infty   P^\mu_{in-1/2}(x) P^\mu_{im-1/2}(x,\pi)  dx$$

$$ = \pi \left[n \sinh(\pi n) \Gamma\left(1/ 2+in-\mu\right) \Gamma \left(1/ 2-in-\mu\right)\right]^{-1} \delta_{n,m},\  n,m \in \mathbb{N},\eqno(2.13)$$
where $\delta_{n,m}$ is the Kronecker delta and integral $(2.13)$ converges absolutely.}

\begin{proof}  In fact,  letting in (1.11) $z=x, \omega=\pi,\ \nu= im-1/2$, we integrate by parts to obtain

$$P^\mu_{\nu}(x, \pi) =  \sqrt{{2\over \pi}}\  {\Gamma(3/2-\mu) (x^2-1)^{-\mu/2} \over m \Gamma(1/2+im-\mu) \Gamma(1/2-im-\mu)} \int_0^\pi  {\sin(mt) \sinh(t) \over (\cosh(t) + x)^{3/2-\mu}}\ dt.\eqno(2.14)$$
Then, plugging the right-hand side of (2.14)  in the left-hand side of (2.13), one can change the order of integration owing to the absolute and uniform convergence since (see the estimate (2.5)) 

$$ \int_0^\pi  \sinh(t)  \int_1^\infty  \left| P^\mu_{in-1/2}(x) \right| { (x^2-1)^{-{\rm Re} \mu/2} \over (x+ \cosh(t))^{3/2-{\rm Re}\mu}} dx dt$$

$$ \le {\pi\over \sqrt 2}  \sinh(\pi)  \int_0^\infty  \left| P^\mu_{in-1/2}(2y+1) \right| { y^{- {\rm Re}\mu/2} \over (y+ 1)^{3/2- {\rm Re}\mu/2}} dy $$

$$\le  {\pi\over \sqrt 2}  \sinh(\pi) C_\mu  \int_0^\infty    { y^{- \gamma - {\rm Re}\mu/2} \over (y+ 1)^{3/2} }dy < \infty\eqno(2.15)$$
for  $\gamma \in  ({\rm Re}\mu /2,\  (1- {\rm Re}\mu) /2) \subset  (- (1+ {\rm Re}\mu)/2,\  1- {\rm Re}\mu /2)$ when $\left| {\rm Re} \mu \right| < 1/2$. 
Therefore
$$\int_1^\infty   P^\mu_{in-1/2}(x)   P^\mu_{im-1/2}(x,\pi)  dx =   \sqrt{{2\over \pi}}\  \frac{\Gamma(3/2-\mu) }{ m\Gamma\left(1/ 2+im-\mu\right) \Gamma \left(1/ 2-im-\mu\right)}$$

$$\times  \int_0^\pi  \sin(mt) \sinh(t)  \int_1^\infty   P^\mu_{in-1/2}(x) { (x^2-1)^{-\mu/2} \over (x+ \cosh(t))^{3/2-\mu}} dx dt.\eqno(2.16)$$
But the inner integral with respect to $x$ on the right-hand side of (2.16) is calculated via formula in [3], Vol. III, Entry 2.17.5.11, and we have

$$\int_1^\infty   P^\mu_{in-1/2}(x) { (x^2-1)^{-\mu/2} \over (x+ \cosh(t))^{3/2-\mu}} dx = {\sqrt {\pi}\over \sqrt 2\  \Gamma(3/2-\mu)\cosh\left(t/2\right) } $$

$$\times {}_2F_1 \left( {1\over 2} -in,\  {1\over 2} +in;\ {1\over 2}; \ \cosh^2\left({t\over 2}\right)\right) - {\sqrt {2\pi}\   n\cosh(\pi n)\over   \sinh(\pi n) \ \Gamma(3/2-\mu) } $$

$$\times  {}_2F_1 \left( 1+in,\  1 -in;\ {3\over 2}; \ \cosh^2\left({t\over 2}\right)\right).\eqno(2.17)$$
Meanwhile,  values of the latter Gauss's hypergeometric functions in (2.17) are calculated in [3], Vol. III, Entries 7.3.1.90 and 7.3.1.94, respectively,

$${}_2F_1 \left( {1\over 2} -in,\  {1\over 2} +in;\ {1\over 2}; \ \cosh^2\left({t\over 2}\right)\right) = { \cos\left( 2in\  \hbox{arcsin} \left(\cosh\left(t/2\right)\right)\right) \over i \sinh(t/2)} ,$$

$${}_2F_1 \left( 1+in,\  1-in;\ {3\over 2}; \ \cosh^2\left({t\over 2}\right)\right) =  - {\sin\left( 2in\  \hbox{arcsin} \left(\cosh\left(t/2\right)\right)\right)\over n \sinh(t)}.$$
Then we derive

$$ \sin\left( 2in\  \hbox{arcsin} \left(\cosh\left(t/2\right)\right)\right)  =  \sin\left( 2n \log\left( i e^{t/2} \right) \right) = \sin\left( nt + \pi n i\right) $$   

$$= \sin(nt) \cosh(\pi n) + i \cos (nt)\sinh(\pi n),$$

$$ \cos \left( 2in\  \hbox{arcsin} \left(\cosh\left(t/2\right)\right)\right)  =  \cos\left( 2n \log\left( i e^{t/2} \right) \right) = \cos\left( nt + \pi n i\right) $$   

$$= \cos(nt) \cosh(\pi n) - i \sin (nt)\sinh(\pi n).$$
Hence, returning to (2.17), it becomes

$$\int_1^\infty   P^\mu_{in-1/2}(x) { (x^2-1)^{-\mu/2} \over (x+ \cosh(t))^{3/2-\mu}} dx = { \sqrt {2\pi}\over \sinh(t) \Gamma(3/2-\mu)}  $$

$$\times \Big[  -i \cos(nt) \cosh(\pi n) -  \sin (nt)\sinh(\pi n)\Big.$$

$$+ \left. { \cosh(\pi n)\over   \sinh(\pi n)} \left[  \sin(nt) \cosh(\pi n) + i \cos (nt)\sinh(\pi n) \right]\right]$$

$$=   { \sqrt {2\pi}\  \sin(nt)  \over  \Gamma(3/2-\mu)  \sinh(t) \sinh(\pi n)}.\eqno(2.18)$$
Therefore, substituting the latter result into (2.16) and making simple calculations, we arrive at (2.13), completing the proof of Lemma 1.  

\end{proof}  

{\bf Corollary 2.} {\it Let  $\left| {\rm Re} \mu \right| < 1/2$. Then under conditions of Theorem $1$ expansion $(1.8)$ takes place}.

\begin{proof}  Since from  (2.5), (2.14) we have the estimate

$$\int_1^\infty   \left| P^\mu_{in-1/2}(x) \right|  $$

$$\times \sum_{m=1}^\infty \left| a_m \Gamma\left({1\over 2} + im -\mu\right) \Gamma\left({1\over 2} - im -\mu\right) P^\mu_{im-1/2}(x, \pi) \right| dx < \infty,$$
the result follows immediately after the interchange of the order of integration and summation on the right-hand side of (1.8) and the use of Lemma 1. 

\end{proof}

Concerning the validity of the expansion (1.9), it is given by the following theorem.

{\bf Theorem 2.} {\it  Let ${\rm Re} \mu < 1/2, \ f$ be a complex-valued function represented by the integral 

$$f(t)=   (t^2-1)^{-\mu/2}  \int_{-\pi}^\pi  {\varphi (u) \ du \over (t+\cosh(u))^{3/2-\mu} },\quad t >1,\eqno(2.19)$$
where  $ \varphi(u) = \psi(u)\sinh(u)$ and $\psi$ is a  $2\pi$-periodic continuously differentiable function, i.e. $\psi \in C^1[-\pi,\pi]$.  Then expansion $(1.9)$ holds for all $x >1$. }

\begin{proof}     Writing (2.19) in the form

$$  f(t)=  { (t^2-1)^{-\mu/2}\over \Gamma(3/2-\mu)}   \int_0^\infty e^{-yt} y^{1/2-\mu} \int_{-\pi}^\pi e^{-y\cosh(u)}  \varphi (u) \  du dy,$$ 
we   treat the integral with respect to $t$ on the right-hand side of (1.9) with the use of formula (2.6) as follows

$$\int_1^\infty  P^\mu_{in-1/2}(t) f(t) dt =  {1\over  \Gamma(3/2-\mu)} \int_1^\infty  P^\mu_{in-1/2}(t) (t^2-1)^{-\mu/2}  \int_0^\infty e^{-yt} y^{1/2-\mu} $$

$$\times \int_{-\pi}^\pi e^{-y\cosh(u)}  \varphi (u) \  du dy dt = {\sqrt 2 \over  \sqrt\pi \Gamma(3/2-\mu)}   \int_0^\infty K_{in}(y) $$

$$\times   \int_{-\pi}^\pi e^{-y\cosh(u)}  \varphi (u) \  du dy,\eqno(2.20)$$
where the interchange of the order of integration is allowed by Fubini's theorem.  Therefore under conditions of the theorem one can invert the latter discrete Kontorovich-Lebedev transform, appealing to Theorem 5 in [6].  As a result we find 

$$\int_{-\pi}^\pi e^{-y\cosh(u)}  \varphi (u) \  du =  {\sqrt 2\over \pi\sqrt\pi}  \Gamma(3/2-\mu) \sum_{n=1}^\infty   n  \sinh(\pi n)  {K_{in} (y,\pi)\over y}$$

$$\times   \int_1^\infty  P^\mu_{in-1/2}(t) f(t) dt,\eqno(2.21) $$
where $K_{in} (y,\pi)$ is defined by (2.9).  Hence, multiplying by $y^{1/2 -\mu} e^{-yt},\ y >0, \ x >1$ both sides of (2.21), we integrate by $y$ over $(0,\infty)$ to obtain via (2.19)

$$ (x^2-1)^{\mu/2} f(x) =   {\sqrt 2\over \pi\sqrt\pi} \int_0^\infty y^{-1/2 -\mu} e^{-yx} $$

$$\times \sum_{n=1}^\infty   n  \sinh(\pi n) K_{in} (y,\pi) \int_1^\infty  P^\mu_{in-1/2}(t) f(t) dt dy.\eqno(2.22)$$
The problem now is to integrate with respect to $y$ under the series sign. To do this we return to (2.20) and,  employing the integral (see [3], Vol. II, Entry  2.16.6.1)

$$ \int_0^\infty  e^{-y\cosh(u)} K_{in}(y)  dy = {\pi \sin( nu) \over \sinh (u) \sinh(\pi n)},\eqno( 2.23)$$
we write the equalities  

$$\int_1^\infty  P^\mu_{in-1/2}(t) f(t) dt = {\sqrt {2\pi} \over \sinh(\pi n)  \Gamma(3/2-\mu)}    \int_{-\pi}^\pi {\varphi (u) \sin(nu) \over  \sinh(u)} du$$

$$=  {\sqrt {2\pi} \over \sinh(\pi n)  \Gamma(3/2-\mu)}    \int_{-\pi}^\pi \psi (u) \sin(nu) du.\eqno(2.24)$$
In the meantime, recalling the incomplete Bessel integral (2.9),  we integrate by parts twice to derive the following representation of the function $K_{in} (y,\pi)$

$$K_{in} (y,\pi) =   { (-1)^{n+1}\over n^2}   \sinh(\pi)\  y e^{-y\cosh(\pi)} $$

$$+   {y\over n^2}  \int_0^\pi  e^{-y\cosh(u)} \left[ \cosh(u)  -  y \sinh^2(u)\right] \cos (n u) du.\eqno(2.25) $$
Moreover, the right-hand side of the latter equality in (2.24) can be written, invoking  the integration by parts. Namely, we find $(\psi (-\pi) = \psi (\pi))$

$$ {\sqrt {2\pi} \over \sinh(\pi n)  \Gamma(3/2-\mu)}    \int_{-\pi}^\pi \psi (u) \sin(nu) du $$

$$=  {\sqrt {2\pi}  \over  n \sinh(\pi n)  \Gamma(3/2-\mu)}  \int_{-\pi}^\pi \psi^\prime (u) \cos(nu) du.\eqno(2.26)$$
Thus, combining with (2.24), (2.25), we recall (2.22) to deduce the estimate of its right-hand side in the form

$$ {\sqrt 2\over \pi\sqrt\pi} \int_0^\infty y^{-1/2 -\mu} e^{-yx} $$

$$\times \sum_{n=1}^\infty   n  \sinh(\pi n) \left| K_{in} (y,\pi) \int_1^\infty  P^\mu_{in-1/2}(t) f(t) dt \right| dy $$

$$\le  {2 \sinh(\pi) \over \pi (x+\cosh(\pi))^{3/2-\mu}}   \int_{-\pi}^\pi \left| \psi^\prime (v) \right| dv \sum_{n=1}^\infty  {1\over n^2}   $$

$$ +  {2 \over \pi  \Gamma(3/2-\mu)}  \int_{-\pi}^\pi \left| \psi^\prime (v) \right| dv \sum_{n=1}^\infty  {1\over n^2}  \int_0^\infty y^{1/2 -\mu} e^{- yx} $$

$$\times   \int_0^\pi  e^{-y\cosh(u)} \left[  \cosh(u)  +  y \sinh^2(u)\right] du  dy$$

$$\le 2 \left[  { \sinh(\pi) \over \pi} +   \cosh(\pi)+ \left( {3\over 2} - \mu\right)  \sinh^2(\pi)\right]
 \int_{-\pi}^\pi \left| \psi^\prime (v) \right| dv \sum_{n=1}^\infty  {1\over n^2} < \infty.\eqno(2.27)$$
 Hence the interchange of the order of integration and summation is permitted, and we get from (2.22), (1.11), (2.9)
 
 $$   f(x) =   {\sqrt 2 (x^2-1)^{- \mu/2}  \over \pi\sqrt\pi} \ \Gamma\left({1\over 2} -\mu\right) $$

$$\times \sum_{n=1}^\infty   n  \sinh(\pi n)  \int_0^\pi {\cos(nu) du \over (x+\cosh(u))^{1/2-\mu} } \int_1^\infty  P^\mu_{in-1/2}(t) f(t) dt$$

$$=  {1\over \pi} \sum_{n=1}^\infty   n  \sinh(\pi n)  \Gamma\left({1\over 2} + in -\mu\right) \Gamma\left({1\over 2} - in -\mu\right) P^\mu_{in-1/2}(x, \pi)$$

$$\times \int_1^\infty  P^\mu_{in-1/2}(t) f(t) dt.$$
Theorem 2 is proved. 

\end{proof} 

{\bf Theorem  3.} {\it  Let $\left|{\rm Re} \mu\right| < 1/2, \ f$ be a complex-valued function represented by the series $(2.10)$, where $ \{a_m\}_{n\in \mathbb{N}} \in l_1$. Then expansion $(1.10)$ holds for all $x > 1$, and the iterated series and integral converge absolutely.}

\begin{proof}  Substituting series (2.10) into the integral in (1.10), we change the order of integration and summation by virtue of the absolute and uniform convergence to derive

$$ \int_1^\infty  P^\mu_{in-1/2}(t,\pi) f(t) dt =   \sum_{m=1}^\infty a_m  \int_1^\infty  P^\mu_{in-1/2}(t,\pi) P^\mu_{im-1/2}(t) dt$$

$$ = \pi  a_n\left[n \sinh(\pi n) \Gamma\left(1/ 2+in-\mu\right) \Gamma \left(1/ 2-in-\mu\right)\right]^{-1}$$
due to the orthogonality (2.13).  Then (1.10) holds.   This proves the theorem.

\end{proof}

\bigskip
\centerline{{\bf Acknowledgments}}
\bigskip

\noindent The work was partially supported by CMUP (UID/MAT/00144/2019),  which is funded by FCT (Portugal) with national (MEC),   European structural funds through the programs FEDER  under the partnership agreement PT2020.

\bigskip
\centerline{{\bf References}}
\bigskip
\baselineskip=12pt
\medskip
\begin{enumerate}

\item[{\bf 1.}\ ]  A. Erd\'elyi,W. Magnus, F. Oberhettinger, and F.G. Tricomi, {\it Higher Transcendental Functions},Vols. I and II, McGraw-Hill, NewYork, London, Toronto, 1953.

\item[{\bf 2.}\ ]  V.A. Fock,  On the representation of an arbitrary function by an integral involving Legendre's functions with a complex index,  {\it C.R. (Doklady) Acad. Sci. URSS (N.S.)} {\bf 39} (1943),  253-256.

\item[{\bf 3.}\ ] A.P. Prudnikov, Yu.A. Brychkov and O.I. Marichev, {\it Integrals and Series}. Vol. I: {\it Elementary Functions}, Vol. II: {\it Special Functions}, Gordon and Breach, New York and London, 1986, Vol. III : {\it More special functions},  Gordon and Breach, New York and London,  1990.

\item[{\bf 4.}\ ]  N.Ya. Vilenkin,  The matrix elements of irreducible unitary representations of a group of Lobachevsky space motions and the generalized Fock-Mehler transformations, {\it Dokl. Akad. Nauk SSSR} (N.S.) {\bf 118},  1958,  219-222.

\item[{\bf 5.}\ ] S. Yakubovich, {\it Index Transforms}, World Scientific Publishing Company, Singapore, New Jersey, London and Hong Kong, 1996.

\item[{\bf 6.}\ ]  S. Yakubovich,   Discrete Kontorovich-Lebedev transforms,  arXiv:1908.01392.

\end{enumerate}

\vspace{5mm}

\noindent S.Yakubovich\\
Department of  Mathematics,\\
Faculty of Sciences,\\
University of Porto,\\
Campo Alegre st., 687\\
4169-007 Porto\\
Portugal\\
E-Mail: syakubov@fc.up.pt\\

\end{document}